\documentclass{amsart}

\usepackage{hyperref}
\usepackage{amsthm,mathrsfs}
\usepackage{pstricks, pst-node}
\usepackage{epsfig}

\newtheorem{thm}{Theorem}[section]
\newtheorem{lem}[thm]{Lemma}

\theoremstyle{definition}

\theoremstyle{remark}

\def\SCH{{\mathrm{SCH}}}
\def\NC#1{{\mathrm{NC}}_2(#1)}
\def\NCp#1{{\mathrm{NC}}_2'(#1)}
\def\CH#1{{\mathrm{CH}}_2(#1)}
\def\CHp#1{{\mathrm{CH}}_2'(#1)}
\def\MOT#1{{\mathrm{MOT}}_3(#1)}
\def\sch.{Schr{\"o}der}
\def\schodd#1{\SCH_{\mathrm{odd}}(#1)}

\def\scheven#1{\SCH_{\mathrm{even}}(#1)}
\def\schevenp#1{\SCH_{\mathrm{even}}'(#1)}
\def\schuh#1{\SCH_{\mathrm{UH}}(#1)}
\def\schuhp#1{\SCH_{\mathrm{UH}}'(#1)}
\def\avoiding#1{P_{12312}(#1)}
\def\avoidingp#1{P_{12312}'(#1)}

\def\w{\mathbf{w}}

\def\bk#1{\{#1\}}

\psset{linewidth=1.5pt,unit=0.7cm, dotsize=5pt}
\psset{subgriddiv=1,gridcolor=lightgray,gridlabels=0pt,gridwidth=0.4pt}

\def\cvput#1[#2]{\pnode(#1,1){#1} \pscircle*(#1,1){.1} \rput(#1,.5){$#2$}}
\def\vput#1{\cvput#1[#1]}
\def\edge#1#2{\ncarc[arcangle=50]{#1}{#2}}
\def\activevertex#1{\pscircle[linecolor=red](#1,1){.2}}

\newcount\sx \newcount\sy \newcount\ex \newcount\ey
\def\CHS#1(#2,#3){ 
\sx=#2 \sy=#3 \ex=#2 \ey=#3
\advance\ex by1 \advance\ey by0
\psline(\sx,\sy)(\ex,\ey)
\rput(\number\sx.5,\number\sy.3){$#1$}
\psdot(\number\sx,\number\sy) \psdot(\number\ex,\number\ey)
}
\def\CDHS#1(#2,#3){ 
\sx=#2 \sy=#3 \ex=#2 \ey=#3
\advance\ex by2 \advance\ey by0
\psline(\sx,\sy)(\ex,\ey)
\rput(\number\sx,\number\sy.3){\rput(1,0){$#1$}}
\psdot(\number\sx,\number\sy) \psdot(\number\ex,\number\ey)
}
\def\CDS#1(#2,#3){ 
\sx=#2 \sy=#3 \ex=#2 \ey=#3
\advance\ex by1 \advance\ey by-1
\psline(\sx,\sy)(\ex,\ey)
\rput(\number\sx.7,\number\ey.7){$#1$}
\psdot(\number\sx,\number\sy) \psdot(\number\ex,\number\ey)
}
\def\HS(#1,#2){ 
\sx=#1 \sy=#2 \ex=#1 \ey=#2
\advance\ex by1 \advance\ey by0
\psline(\sx,\sy)(\ex,\ey)
\psdot(\number\sx,\number\sy) \psdot(\number\ex,\number\ey)
}
\def\US(#1,#2){ 
\sx=#1 \sy=#2 \ex=#1 \ey=#2
\advance\ex by1 \advance\ey by1
\psline(\sx,\sy)(\ex,\ey)
\psdot(\number\sx,\number\sy) \psdot(\number\ex,\number\ey)
}
\def\CUS#1(#2,#3){ 
\sx=#2 \sy=#3 \ex=#2 \ey=#3
\advance\ex by1 \advance\ey by1
\psline(\sx,\sy)(\ex,\ey)
\rput(\number\sx.3,\number\sy.7){$#1$}
\psdot(\number\sx,\number\sy) \psdot(\number\ex,\number\ey)
}
\def\DS(#1,#2){ 
\sx=#1 \sy=#2 \ex=#1 \ey=#2
\advance\ex by1 \advance\ey by-1
\psline(\sx,\sy)(\ex,\ey)
\psdot(\number\sx,\number\sy) \psdot(\number\ex,\number\ey)
}
\def\DHS(#1,#2){ 
\sx=#1 \sy=#2 \ex=#1 \ey=#2
\advance\ex by2 \advance\ey by0
\psline(\sx,\sy)(\ex,\ey)
\psdot(\number\sx,\number\sy) \psdot(\number\ex,\number\ey)
}
\def\CURL(#1,#2)#3{ 
\rput(#1,#2){
\rput(1.5,1.3){#3}
\rput(1.5,0.3){$\cdots$}
\rput(1.5,0.75){\psscaleboxto(2.5,.5){\rotateright{\{}}}
}}

\begin{document}

\title{Bijections on two variations of noncrossing partitions}

\author{Jang Soo Kim}
\email{jskim@kaist.ac.kr}

\date{\today}

\begin{abstract}
We find bijections on 2-distant noncrossing partitions, 12312-avoiding
partitions, 3-Motzkin paths, UH-free Schr{\"o}der paths and
Schr{\"o}der paths without peaks at even height.  We also give a
direct bijection between 2-distant noncrossing partitions and
12312-avoiding partitions.
\end{abstract}

\maketitle

\section{Introduction}

Noncrossing partitions were first introduced by Kreweras \cite{Kreweras1972} in
1972. Recently, they have received great attention, and have been generalized in
many different ways; for instance, see \cite{Armstrong, Chen2007,Drake2008,
  Kasraoui2006,Mansour2007} and the references therein.  In this paper we
consider two variations of noncrossing partitions: $k$-distant noncrossing
partitions and $12\cdots r12$-avoiding partitions introduced by Drake and Kim
\cite{Drake2008}, and Mansour and Severini \cite{Mansour2007} respectively,
where they reduce to noncrossing partitions when $k=1$ and $r=2$.
 
\begin{figure}
  \begin{center}
\psset{linewidth=1pt}
\begin{pspicture}(1,.5)(9,2)
\vput{1} \vput{2} \vput{3} \vput{4} \vput{5} \vput{6} \vput{7}
\vput{8} \vput{9} \edge{1}{4} \edge{4}{8} \edge{2}{5} \edge{5}{9}
\edge{6}{7}
\end{pspicture}
  \end{center}
  \caption{The standard representation of
    $(\bk{1,4,8},\bk{2,5,9},\bk{3},\bk{6,7})$.}
\label{fig:diagram}
\end{figure}
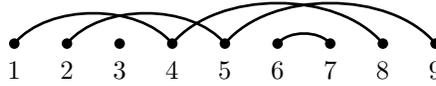

A (set) {\em partition} of $[n]=\{1,2,\ldots,n\}$ is a collection of mutually
disjoint nonempty subsets, called \emph{blocks}, of $[n]$ whose union is $[n]$.
We will write a partition as a sequence of blocks $(B_1,B_2,\ldots,B_k)$ such
that $\min(B_1)<\min(B_2)<\cdots<\min(B_k)$.  An \emph{edge} of a partition is a
pair $(i,j)$ of integers contained in the same block that does not contain any
integer $t$ with $i<t<j$. The \emph{standard representation} of a partition
$\pi$ of $[n]$ is the diagram having $n$ vertices labeled with $1,2,\dots,n$,
where $i$ and $j$ are connected by an arc if $(i,j)$ is an edge of $\pi$; see
Figure~\ref{fig:diagram}.  A \emph{noncrossing partition} is a partition without
any two crossing edges, i.e. $(i_1,j_1)$ and $(i_2,j_2)$ such that
$i_1<i_2<j_1<j_2$. It is well known the number of noncrossing partitions of
$[n]$ is the Catalan number $\frac{1}{n+1}\binom{2n}{n}$.

For a positive integer $k$, a \emph{$k$-distant noncrossing partition} is a
partition without any two edges $(i_1, j_1)$ and $(i_2, j_2)$ satisfying
$i_1<i_2 < j_1<j_2$ and $j_1 - i_2 \geq k$. Note that $1$-distant noncrossing
partitions are just noncrossing partitions. We denote by $\mathrm{NC}_k(n)$ the
set of $k$-distant noncrossing partitions of $[n]$.  Drake and Kim
\cite{Drake2008} found the following generating function for the number of
$2$-distant noncrossing partitions:
\begin{equation}\label{eq:gf}
\sum_{n\geq0} \# \NC{n} x^n = \frac{3-3x-\sqrt{1-6x+5x^2}}{2(1-x)}.
\end{equation}

The \emph{canonical word} of a partition $\pi=(B_1,B_2,\ldots,B_k)$ is the word
$a_1a_2\cdots a_n$, where $a_i=j$ if $i\in B_j$. For instance, the canonical
word of the partition in Figure~\ref{fig:diagram} is $123124412$. In the
literature canonical words are also called restricted growth functions. For a
word $\tau$, a partition is called \emph{$\tau$-avoiding} if its canonical word
does not contain a subword which is order-isomorphic to $\tau$. It is easy to
see that a partition is noncrossing if and only if it is $1212$-avoiding. We
denote by $P_\tau(n)$ the set of $\tau$-avoiding partitions of $[n]$.

Using the kernel method, Mansour and Severini \cite{Mansour2007} found the
generating function for the number of $12\cdots r12$-avoiding partitions of
$[n]$. Interestingly, as a special case of their result, the generating function
for the number of $12312$-avoiding partitions of $[n]$ is the same as
\eqref{eq:gf}, which implies $\#\NC{n}=\#\avoiding{n}$.  Moreover, this number
also counts several kinds of lattice paths.  The main purpose of this paper is
to find bijections between $\NC n$ and $\avoiding n$ together with some lattices
paths described below.

A \emph{lattice path of length} $n$ is a sequence of points in $\mathbb{N}\times
\mathbb{N}$ starting at $(0,0)$ and ending at $(n,0)$. For a lattice path
$L=((x_0,y_0),(x_1,y_1),\ldots,(x_k,y_k))$, each $S_i=(x_i-x_{i-1},
y_i-y_{i-1})$ is called a \emph{step} of $L$.  The \emph{height} of the step
$S_i$ is defined to be $y_{i-1}$. Sometimes we will identify a lattice path $L$
with the word $S_1S_2\ldots S_k$ of its steps. Note that the number of steps is
not necessarily equal to the length of the lattice path.

Let $U$, $D$ and $H$ denote an \emph{up step}, a \emph{down step} and
a \emph{horizontal step} respectively, i.e., $U=(1,1)$, $D=(1,-1)$ and
$H=(1,0)$.

A \emph{\sch. path} is a lattice path consisting of steps $U$, $D$ and
$H^2=HH=(2,0)$.  Let $L=S_1 S_2\cdots S_k$ be a \sch. path.  A \emph{UH-pair} of
$L$ is a pair $(S_i,S_{i+1})$ of consecutive steps such that $S_i=U$ and
$S_{i+1}=H^2$.  We say that $L$ is \emph{UH-free} if it does not have a UH-pair.
A \emph{peak} of $L$ is a pair $(S_i,S_{i+1})$ of consecutive steps such that
$S_i=U$ and $S_{i+1}=D$. The \emph{height} of a peak $(S_i,S_{i+1})$ is the
height of $S_{i+1}=D$.  We denote by $\schuh{n}$ the set of UH-free \sch. paths
of length $2n$, and by $\scheven{n}$ (resp. $\schodd{n})$ the set of \sch. paths
of length $2n$ which have no peaks of even (resp. odd) height.

A \emph{labeled step} is a step together with an integer label.  Let $D_i$
(resp. $H_i$) denote a labeled down step (resp. a labeled horizontal step) with
label $i$.  We denote by $\CH{n}$ the set of lattice paths $L=S_1S_2\cdots S_n$
of length $n$ consisting of $U$, $D_1$, $D_2$, $H_0$, $H_1$ and $H_2$ such that
\begin{itemize}
\item if $S_i=H_{\ell}$ or $S_i=D_{\ell}$, then $S_i$ is of height at
  least $\ell$,
\item if $S_i=H_2$ or $S_i=D_2$, then $i\geq2$ and
$S_{i-1}\in\{U,H_1,H_2\}$.
\end{itemize}

A \emph{3-Motzkin path} is a lattice path consisting of $U$, $D$, $H_0$, $H_1$
and $H_2$.  We denote by $\MOT{n}$ the set of 3-Motzkin paths of length $n$.

Drake and Kim \cite{Drake2008} showed that the well known bijection $\psi$
between partitions and Charlier diagrams, see \cite{Flajolet1980,Kasraoui2006},
yields a bijection $\psi:\NC{n} \rightarrow \CH{n}$.  Yan \cite{Yan2009} found a
bijection $\phi:\schuh{n-1}\rightarrow \avoiding{n}$ and a bijection between
$\schuh{n}$ and $\scheven{n}$.  Thus all of $\NC{n}$, $\CH{n}$, $\scheven{n-1}$,
$\schuh{n-1}$ and $\avoiding{n}$ have the same cardinality, which is counted by
sequence A007317 from \cite{Sloane}.  In order to find bijections between these
objects, we introduce the following sets:
\begin{itemize}
\item $\NCp{n}=\{\pi\in \NC{n} : \mbox{$n$ is not a singleton}\}$
\item $\CHp{n}=\{L\in \CH{n} : \mbox{the last step of $L$ is $D_1$}\}$
\item $\schevenp{n}=\{L\in \scheven{n} : \mbox{the first step of $L$
  is $U$}\}$
\item $\schuhp{n}=\{L\in \schuh{n} : \mbox{the first step of $L$ is
  $U$}\}$
\item $\avoidingp{n}=\{\pi\in \avoiding{n} : \mbox{$1$ and $2$ are not
in the same block}\}$
\end{itemize}

Note that we can identify $\pi\in \NC n$ with $\pi'\in \NCp k$, where $k$ is the
integer such that $j$ is a singleton for all $j\in \{k+1,k+2,\dots,n\}$ and $k$
is not a singleton in $\pi$, and $\pi'$ is the partition obtained from $\pi$ by
deleting integers greater than $k$. We can also identify $\pi\in \avoiding n$
with $\overline{\pi}\in \avoidingp k$, where $k$ is the integer such that the
number of consecutive $1$'s at the beginning of the canonical word of $\pi$ is
$n-k+1$, and $\overline{\pi}$ is the partition whose canonical word is obtained
from that of $\pi$ by deleting the first $n-k$ $1$'s. Thus any bijection between
$\NCp n$ and $\avoidingp n$ naturally induces a bijection between $\NC n$ and
$\avoiding n$. Similarly, any bijection between $A'(n)$ and $B'(n)$ naturally
induces a bijection between $A(n)$ and $B(n)$ where $A$ and $B$ are any two of
$\mathrm{NC}_2$, $\mathrm{CH}_2$, $\mathrm{SCH_{even}}$, $\mathrm{SCH_{UH}}$,
and $P_{12312}$. Thus in order to find a bijection between $\NC n$ and
$\avoiding n$, it is enough to find a bijection between $\NCp n$ and $\avoidingp
n$. 

In this paper we find bijections between these objects.  For the overview of our
bijections see Figure~\ref{fig:bijections}, where $\psi$ is the known bijection
between partitions and Charlier diagrams \cite{Flajolet1980,Kasraoui2006}, and
$\phi$ is Yan's bijection \cite{Yan2009}. We note that our bijection $g$ in
Figure~\ref{fig:bijections} is also discovered by Shapiro and Wang
\cite{Shapiro2009}. We also provide a direct bijection between $\NC{n}$ and
$\avoiding{n}$ in Section~\ref{sec:direct}.

\begin{figure}
\[
\psset{linewidth=1pt,arrows=->,nodesep=4pt,colsep=1cm, rowsep=1cm,
shortput=nab}
\psmatrix
\NCp{n} & & \schevenp{n-1}& \avoidingp{n}\\
\CHp{n} & \MOT{n-2} & \schodd{n-1} & \schuhp{n-1}
\endpsmatrix
\ncline{1,1}{2,1}^{\psi}
\ncline{2,1}{2,2}^{f}
\ncline{2,2}{2,3}^{g}
\ncline{2,3}{2,4}^{h}
\ncline{2,3}{1,3}_{\iota}
\ncline{2,4}{1,4}_{\phi}
\]
\caption{Main bijections for $n\geq2$.} \label{fig:bijections}
\end{figure}

\section{Bijections}\label{sec:bijections}

In this section we find the bijections $f,g,h$, and $\iota$ in
Figure~\ref{fig:bijections}.

\subsection{The bijection $f:\CHp{n}\rightarrow \MOT{n-2}$}

Recall that $\CHp{n}$ is the set of lattice paths $L=S_1S_2\cdots S_n$
of length $n$ consisting of $U, D_1, D_2, H_0, H_1$ and $H_2$ such
that
\begin{itemize}
\item if $S_i=H_{\ell}$ or $S_i=D_{\ell}$, then $S_i$ is of height at
  least $\ell$,
\item if $S_i=H_2$ or $S_i=D_2$, then $i\geq2$ and
$S_{i-1}\in\{U,H_1,H_2\}$,
\item $S_n=D_1$.
\end{itemize}
The second condition above is equivalent to the condition that the
lattice path consists of the following combined steps for any $k\geq
0$:
\begin{equation}\label{eq:steps}
UH_2^k,UH_2^kD_2, H_1H_2^k, H_1H_2^kD_2, H_0,D_1.
\end{equation}

Let $A(n)$ denote the set of lattice paths of length $n$ consisting of
the combined steps in \eqref{eq:steps} such that $H_2$ does not touch
the $x$-axis.  Let $B(n)$ denote the set of 3-Motzkin paths of
length $n$ such that each $H_2$ touching the $x$-axis must occur after
$D$, $H_0$ or $H_2$.

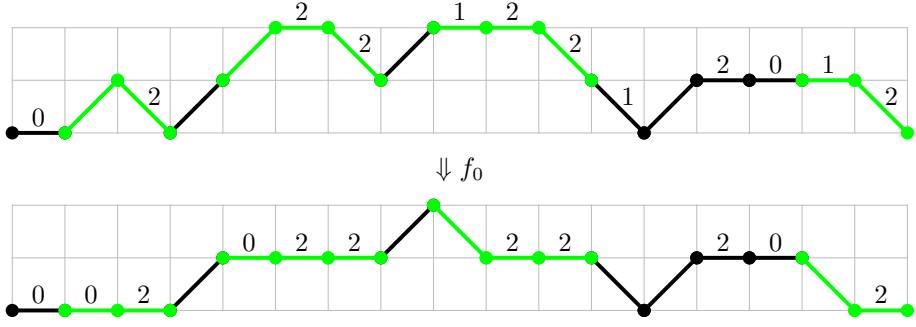
\begin{figure}
\begin{center}
\begin{pspicture}(0,0)(17,2.2)
\psgrid(0,0)(17,2) \CHS{0}(0,0) \US(3,0) \US(7,1) \US(12,0)
\CHS{2}(13,1) \CHS{0}(14,1) \CDS{1}(11,1) \psset{linecolor=green}
\US(1,0) \CDS{2}(2,1) \US(4,1) \CHS{2}(5,2) \CDS{2}(6,2) \CHS{1}(8,2)
\CHS{2}(9,2) \CDS{2}(10,2) \CHS{1}(15,1) \CDS{2}(16,1)
\end{pspicture}
\[\Downarrow f_0\]
\begin{pspicture}(0,0)(17,2.2)
\psgrid(0,0)(17,2) \CHS{0}(0,0) \US(3,0) \US(7,1) \DS(11,1) \US(12,0)
\CHS{2}(13,1) \CHS{0}(14,1) \psset{linecolor=green} \CHS{0}(1,0)
\CHS{2}(2,0) \CHS{0}(4,1) \CHS{2}(5,1) \CHS{2}(6,1) \DS(8,2)
\CHS{2}(9,1) \CHS{2}(10,1) \DS(15,1) \CHS{2}(16,0)
\end{pspicture}
\end{center}
\caption{An example of $f_0$.}
\label{fig:f0}
\end{figure}

We define $f_0:A(n)\rightarrow B(n)$ as follows. Let $L\in A(n)$.
Then $f_0(L)$ is defined to be the lattice path obtained from $L$ by
changing $UH_2^kD_2$ to $H_0H_2^{k+1}$, $H_1H_2^kD_2$ to $DH_2^{k+1}$
and $D_1$ to $D$.  It is easy to see that $f_0(L)\in B$ and $f_0$ is
invertible. See Figure~\ref{fig:f0}.

\begin{figure}
\begin{center}
\scalebox{0.7}{\begin{pspicture}(0,0)(25,3)
    \CURL(0,0){$k_1$} \CHS{0}(0,0) \HS(1,0) \CHS{0}(2,0) \US(3,0)
    \pswedge(5.5,1){1.5}{0}{180}
    \rput(5.5,1.7){$L_1$}\DS(7,1) \CURL(8,0){$k_2$} \CHS{0}(8,0)
    \HS(9,0) \CHS{0}(10,0) \US(11,0)
    \pswedge(13.5,1){1.5}{0}{180}
    \rput(13.5,1.7){$L_2$}\DS(15,1) \CURL(17,0){$k_r$} \CHS{0}(17,0)
    \HS(18,0) \CHS{0}(19,0) \US(20,0) \rput(16.5,0.5){$\cdots$}
    \pswedge(22.5,1){1.5}{0}{180}
    \rput(22.5,1.7){$L_r$}\DS(24,1)
\end{pspicture}}
\[\Downarrow f\]
\scalebox{0.7}{\begin{pspicture}(0,0)(23,2)
\CURL(0,0){$k_1$} \CHS{2}(0,0) \HS(1,0)
\CHS{2}(2,0) \pswedge(4.5,0){1.5}{0}{180}
\rput(4.5,0.7){$f_0(L_1)$}
\CHS{1}(6,0) \CHS{2}(7,0) \CURL(8,0){$k_2$}
\CHS{2}(8,0) \HS(9,0) \CHS{2}(10,0)
\pswedge(12.5,0){1.5}{0}{180}
\rput(12.5,0.7){$f_0(L_2)$}
\CHS{1}(15,0) \CHS{2}(16,0) \CURL(17,0){$k_r$}
\CHS{2}(17,0) \HS(18,0) \CHS{2}(19,0)
\pswedge(21.5,0){1.5}{0}{180}
\rput(21.5,0.7){$f_0(L_r)$}
\rput(14.5,0.5){$\cdots$}
\end{pspicture}}
\end{center}
  \caption{Definition of $f$.}\label{fig:f}
\end{figure}
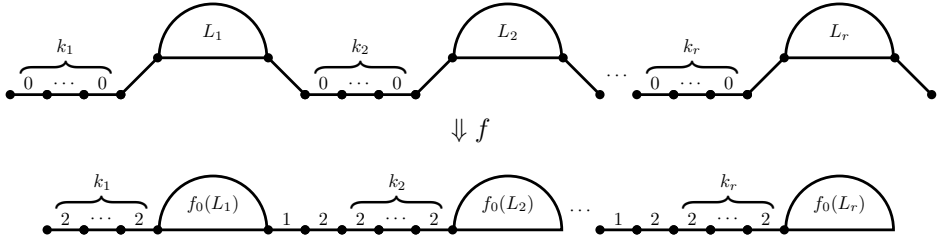

Now we define $f:\CHp{n}\rightarrow \MOT{n-2}$ as follows.  Let $L\in
\CHp{n}$. Then $L$ is decomposed uniquely as
\[
H_0^{k_1} (U L_1 D_1) H_0^{k_2} (U L_2 D_1) \cdots H_0^{k_r} (U L_r D_1),
\]
 where $L_i\in A(n_i)$ for some $k_i,n_i\geq0$ and $r\geq1$.
Then define $f(L)$ to be
\[H_2^{k_1} f_0(L_1) (H_1H_2^{k_2+1} f_0(L_2)) (H_1 H_2^{k_3+1}
f_0(L_3)) \cdots (H_1 H_2^{k_r+1} f_0(L_r)).\] See Figure~\ref{fig:f}.

\begin{thm}
The map $f:\CHp{n}\rightarrow \MOT{n-2}$ is a bijection.
\end{thm}
\begin{proof}
Each $L\in\MOT{n-2}$ is uniquely decomposed as
\[H_2^{k_1} L_1 (H_1 H_2^{k_2+1} L_2) (H_1 H_2^{k_3+1} L_3) \cdots
(H_1 H_2^{k_r+1} L_r),\] where $L_i\in B(n_i)$ for some
$k_i,n_i\geq 0$ and $r\geq 1$. Thus we have the inverse $f^{-1}(L)$
which is decomposed as
\[H_0^{k_1} (U f_0^{-1}(L_1) D_1) H_0^{k_2} (U f_0^{-1}(L_2) D_1)
\cdots H_0^{k_r} (U f_0^{-1}(L_r) D_1).\]
\end{proof}

\subsection{The bijection $g:\MOT{n}\rightarrow \schodd{n+1}$}
\begin{figure}
\def\oddpeak(#1,#2){
\pscircle[linewidth=1.5pt](#1,#2){6pt}}
  \begin{center}
\begin{pspicture}(0,0)(13,2)
\psgrid(0,0)(13,2) \CHS{0}(0,0) \US(1,0) \DS(2,1) \US(3,0)
\CHS{2}(4,1) \DS(5,1) \CHS{1}(6,0) \US(7,0) \US(8,1) \DS(9,2)
\DS(10,1) \CHS{1}(11,0) \CHS{1}(12,0)
\end{pspicture}
\[\downarrow\]
\scalebox{0.5}{%
\psset{linewidth=2pt}
\begin{pspicture}(0,0)(28,6)
\psgrid[gridwidth=0.8pt](0,0)(28,6) \US(0,0) \DHS(1,1) \US(3,1)
\US(4,2) \DS(5,3) \DS(6,2) \US(7,1) \US(8,2) \US(9,3) \DS(10,4)
\DS(11,3) \DS(12,2) \DS(13,1) \US(14,0) \US(15,1) \US(16,2) \US(17,3)
\US(18,4) \DS(19,5) \DS(20,4) \DS(21,3) \DS(22,2) \DS(23,1) \US(24,0)
\DS(25,1) \US(26,0) \DS(27,1) \oddpeak(5,3) \oddpeak(19,5)
\oddpeak(25,1) \oddpeak(27,1)
\end{pspicture}}
\[\downarrow\] 
\scalebox{0.5}{%
\psset{linewidth=2pt}
\begin{pspicture}(0,0)(28,6)
\psgrid[gridwidth=0.8pt](0,0)(28,6) \US(0,0) \DHS(1,1) \US(3,1)
\DS(6,2) \US(7,1) \US(8,2) \US(9,3) \DS(10,4) \DS(11,3)
\DS(12,2) \DS(13,1) \US(14,0) \US(15,1) \US(16,2) \US(17,3)
\DS(20,4) \DS(21,3) \DS(22,2) \DS(23,1) \psset{linecolor=blue, linestyle=dashed}
\DHS(4,2) \DHS(18,4) \DHS(24,0) \DHS(26,0)
\end{pspicture}}
  \end{center}
  \caption{An example of $g$. Odd peaks are circled. The horizontal steps of
    even height are dashed and colored blue.}\label{fig:g}
\end{figure}
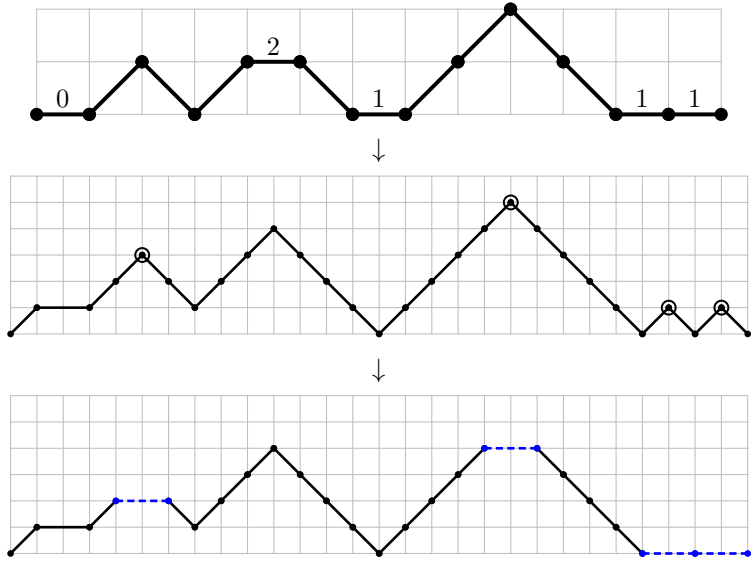

We define $g:\MOT{n}\rightarrow \schodd{n+1}$ as follows.  Let
$L\in\MOT{n}$. Then $g(L)$ is the lattice path obtained from $L$ by
doing the following.
\begin{enumerate}
\item Change $U$ to $UU$, $D$ to $DD$, $H_0$ to $H^2$, $H_1$ to $DU$,
  and $H_2$ to $UD$.
\item Add $U$ at the beginning and $D$ at the end.
\item Change all the consecutive steps $UD$ which form a peak of odd
  height to $H^2$.
\end{enumerate}
See Figure~\ref{fig:g} for an example of $g$.

\begin{thm}
The map $g:\MOT{n}\rightarrow \schodd{n+1}$ is a bijection.
\end{thm}
\begin{proof}
  Clearly the first and the second steps in the construction of $g$
  are invertible. The third step is also invertible because every step
  $H^2$ of even height always comes from a peak of odd height. Thus
  $g$ is invertible.
\end{proof}

\subsection{The bijection $h:\schodd{n}\rightarrow \schuhp{n}$}
Let $L=S_1S_2\cdots S_k$ be a \sch. path.  For any up step $S_i=U$ of $L$, there
is a unique down step $S_j=D$ such that $i<j$ and $S_{i+1}S_{i+2}\cdots S_{j-1}$
is a (possibly empty) lattice path. We call such $S_j$ \emph{the down step
  corresponding} to $S_i$.  We also call $S_i$ \emph{the up step corresponding}
to $S_j$.

For a UH-pair $(S_{i},S_{i+1})$, i.e.  $S_i=U$ and $S_{i+1}=H^2$, 
we define the function $\xi$ as follows.
\[\xi(S_i,S_{i+1})=\left\{\begin{array}{ll} i, & \mbox{if $S_{i+1}$ is
  of even height;}\\ j, & \mbox{if $S_{i+1}$ is of odd height,}
  \end{array}\right.\]
where $j$ is the integer such that $S_j$ is the down step corresponding to
$S_i$.  If $L$ is not UH-free, we define the \emph{$\xi$-maximal} UH-pair of $L$
to be the UH-pair $(S_i,S_{i+1})$ with the largest $\xi$ value.

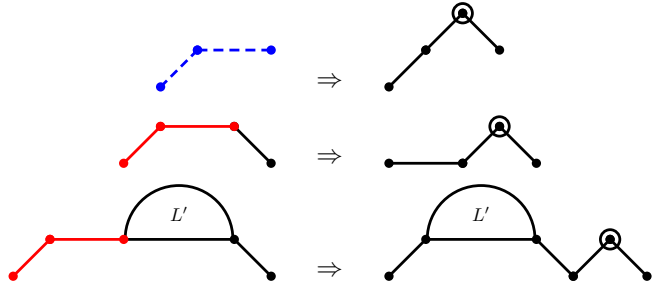
\begin{figure}
  \begin{center}
    \scalebox{0.7}{\begin{pspicture}(0,0)(4,2)
      \psset{linecolor=blue, linestyle=dashed}
      \US(0,0) \DHS(1,1) 
    \end{pspicture}}
    $\Rightarrow$
    \scalebox{0.7}{\begin{pspicture}(-1,0)(3,2)
      \US(0,0) \US(1,1) \DS(2,2) \pscircle(2,2){6pt}
    \end{pspicture}}

    \scalebox{0.7}{\begin{pspicture}(0,0)(5,2)
      \DS(3,1) \psset{linecolor=red} \US(0,0) \DHS(1,1)
    \end{pspicture}}
    $\Rightarrow$
    \scalebox{0.7}{\begin{pspicture}(-1,0)(4,1)
      \DHS(0,0) \US(2,0) \DS(3,1) \pscircle(3,1){6pt}
    \end{pspicture}}

    \scalebox{0.7}{\begin{pspicture}(0,0)(8,3)
      \pswedge(4.5,1){1.5}{0}{180} \rput(4.5,1.7){$L'$}\DS(6,1)
      \psset{linecolor=red} \US(0,0) \DHS(1,1) 
    \end{pspicture}}
    $\Rightarrow$
    \scalebox{0.7}{\begin{pspicture}(-1,0)(7,3)
      \US(0,0) \pswedge(2.5,1){1.5}{0}{180} \rput(2.5,1.7){$L'$}\DS(4,1) \US(5,0) \DS(6,1) 
      \pscircle(6,1){6pt}
    \end{pspicture}}
  \end{center}
  \caption{The essence of $h_0$. Red (resp.~Dashed blue) color is for UH-pairs
    whose horizontal step is of odd (resp.~even) height. Odd peaks are
    circled. The lattice path $L'$ is not empty.}\label{fig:uh}
\end{figure}

Now let $L=S_1S_2\cdots S_k\in \schodd n$. If $L$ is not UH-free, we define
$h_0(L)$ as follows.  Suppose $(S_{i},S_{i+1})$ is the $\xi$-maximal UH-pair of
$L$, and $S_j$ is the down step corresponding to $S_i$.
\begin{enumerate}
\item If $S_{i+1}$ is of even height, then $h_0(L)$ is the lattice
  path obtained from $L$ by replacing $S_{i}S_{i+1}$ with $UUD$.
\item If $S_{i+1}$ is of odd height, then let $L'=S_{i+2}S_{i+3}\cdots
  S_{j-1}$.
\begin{enumerate}
\item If $L'$ is empty, i.e., $j=i+2$, then $h_0(L)$ is the lattice path
obtained from $L$ by replacing $S_{i}S_{i+1}S_{i+2}$ with $H^2UD$.
\item If $L'$ is not empty, then $h_0(L)$ is the lattice path obtained
  from $L$ by replacing $S_{i}S_{i+1} \cdots S_{j}$ with $UL'DUD$.
\end{enumerate}
\end{enumerate}
See Figure~\ref{fig:uh}.

\begin{figure}
  \begin{center}
\psset{unit=0.5cm, dotsize=3pt, linewidth=1pt, gridwidth=0.5pt}
\def\oddpeak(#1,#2){
\pscircle[linewidth=1pt](#1,#2){4pt}}
\def\odduh(#1,#2)#3{
\rput(#1.5,#2.5){
\psline[linestyle=dashed,linewidth=1pt,linecolor=gray]{->}(0,0)(#3,0)}}
\begin{pspicture}(-2,0)(22,3)
\psgrid(-2,0)(22,3) \US(0,0) \DHS(1,1) \US(3,1) \DS(4,2) \DHS(-2,0)
\US(8,2) \DHS(9,3) \DS(11,3) \DS(12,2) \DS(13,1)
\US(14,0) \DHS(15,1) \DHS(17,1) \DS(19,1) \DHS(20,0)
   {\psset{linecolor=blue,linestyle=dashed}\US(5,1) \DHS(6,2)} {\psset{linecolor=red}
     \US(14,0) \DHS(15,1) \US(0,0) \DHS(1,1) \US(8,2) \DHS(9,3)}
   \odduh(14,0)5 \odduh(8,2)3 \odduh(0,0){13}
\end{pspicture}
\[\Downarrow h_0\]
\begin{pspicture}(-2,0)(22,3)
\psgrid(-2,0)(22,3) \US(0,0) \DHS(1,1) \US(3,1) \DS(4,2) \DHS(-2,0)
\US(8,2) \DHS(9,3) \DS(11,3) \DS(12,2) \DS(13,1)
\US(14,0) \DHS(15,1) \DS(17,1) \US(18,0) \DS(19,1) \DHS(20,0)
   {\psset{linecolor=blue,linestyle=dashed}\US(5,1) \DHS(6,2)} {\psset{linecolor=red}
     \US(14,0) \DHS(15,1) \US(0,0) \DHS(1,1) \US(8,2) \DHS(9,3)}
   \oddpeak(19,1) \odduh(14,0)3 \odduh(8,2)3 \odduh(0,0){13}
\end{pspicture}
\[\Downarrow h_0\]
\begin{pspicture}(-2,0)(22,3)
\psgrid(-2,0)(22,3) \US(0,0) \DHS(1,1) \US(3,1) \DS(4,2) \DHS(-2,0)
\US(8,2) \DHS(9,3) \DS(11,3) \DS(12,2) \DS(13,1)
\DHS(14,0) \US(16,0) \DS(17,1) \US(18,0) \DS(19,1) \DHS(20,0)
    {\psset{linecolor=blue,linestyle=dashed}\US(5,1) \DHS(6,2)} {\psset{linecolor=red}
      \US(8,2) \DHS(9,3) \US(0,0) \DHS(1,1)} \oddpeak(17,1)
    \oddpeak(19,1) \odduh(8,2)3 \odduh(0,0){13}
\end{pspicture}
\[\Downarrow h_0\]
\begin{pspicture}(-2,0)(22,3) \psgrid(-2,0)(22,3) \US(0,0) \US(1,1) \DS(2,2)
  \DHS(-2,0) \US(6,2) \DHS(7,3) \DS(9,3) \DS(10,2) \DS(11,1) \US(12,0) \DS(13,1)
  \DHS(14,0) \US(16,0) \DS(17,1) \US(18,0) \DS(19,1)\DHS(20,0)
  {\psset{linecolor=blue,linestyle=dashed} \US(3,1) \DHS(4,2)}
  {\psset{linecolor=red} \US(6,2) \DHS(7,3)}
  \oddpeak(13,1)\oddpeak(17,1)\oddpeak(19,1) \odduh(6,2)3
\end{pspicture}
\[\Downarrow h_0\]
\begin{pspicture}(-2,0)(22,3)
\psgrid(-2,0)(22,3) \US(0,0) \US(1,1) \DS(2,2) \DHS(-2,0) 
\DHS(6,2) \US(8,2) \DS(9,3) \DS(10,2) \DS(11,1) \US(12,0)
\DS(13,1) \DHS(14,0) \US(16,0) \DS(17,1) \US(18,0) \DS(19,1)\DHS(20,0)
   {\psset{linecolor=blue,linestyle=dashed} \US(3,1) \DHS(4,2)} \oddpeak(9,3)
   \oddpeak(13,1) \oddpeak(17,1) \oddpeak(19,1)
\end{pspicture}
\[\Downarrow h_0\]
\begin{pspicture}(-2,0)(22,3)
\psgrid(-2,0)(22,3) \US(0,0) \US(1,1) \DS(2,2) \DHS(-2,0) \US(3,1)
\US(4,2) \DS(5,3) \DHS(6,2) \US(8,2) \DS(9,3) \DS(10,2) \DS(11,1)
\US(12,0) \DS(13,1) \DHS(14,0) \US(16,0) \DS(17,1) \US(18,0)
\DS(19,1)\DHS(20,0) \oddpeak(5,3) \oddpeak(9,3) \oddpeak(13,1)
\oddpeak(17,1) \oddpeak(19,1)
\end{pspicture}
\[\Downarrow \mbox{Changing the first $H^2$}\]
\begin{pspicture}(-2,0)(22,3)
\psgrid(-2,0)(22,3) \US(0,0) \US(1,1) \DS(2,2) \US(-2,0)\DS(-1,1)
\US(3,1) \US(4,2) \DS(5,3) \DHS(6,2) \US(8,2) \DS(9,3) \DS(10,2)
\DS(11,1) \US(12,0) \DS(13,1) \DHS(14,0) \US(16,0) \DS(17,1) \US(18,0)
\DS(19,1)\DHS(20,0) \oddpeak(-1,1) \oddpeak(5,3) \oddpeak(9,3)
\oddpeak(13,1) \oddpeak(17,1) \oddpeak(19,1)
\end{pspicture}
  \end{center}
  \caption{An example of $h$.  Red (resp. Dashed blue) color is for UH-pairs
    whose horizontal step is of odd (resp.~even) height. Odd peaks are
    circled. Dashed arrows indicate the down steps corresponding to the up
    steps.}\label{fig:h}
\end{figure}
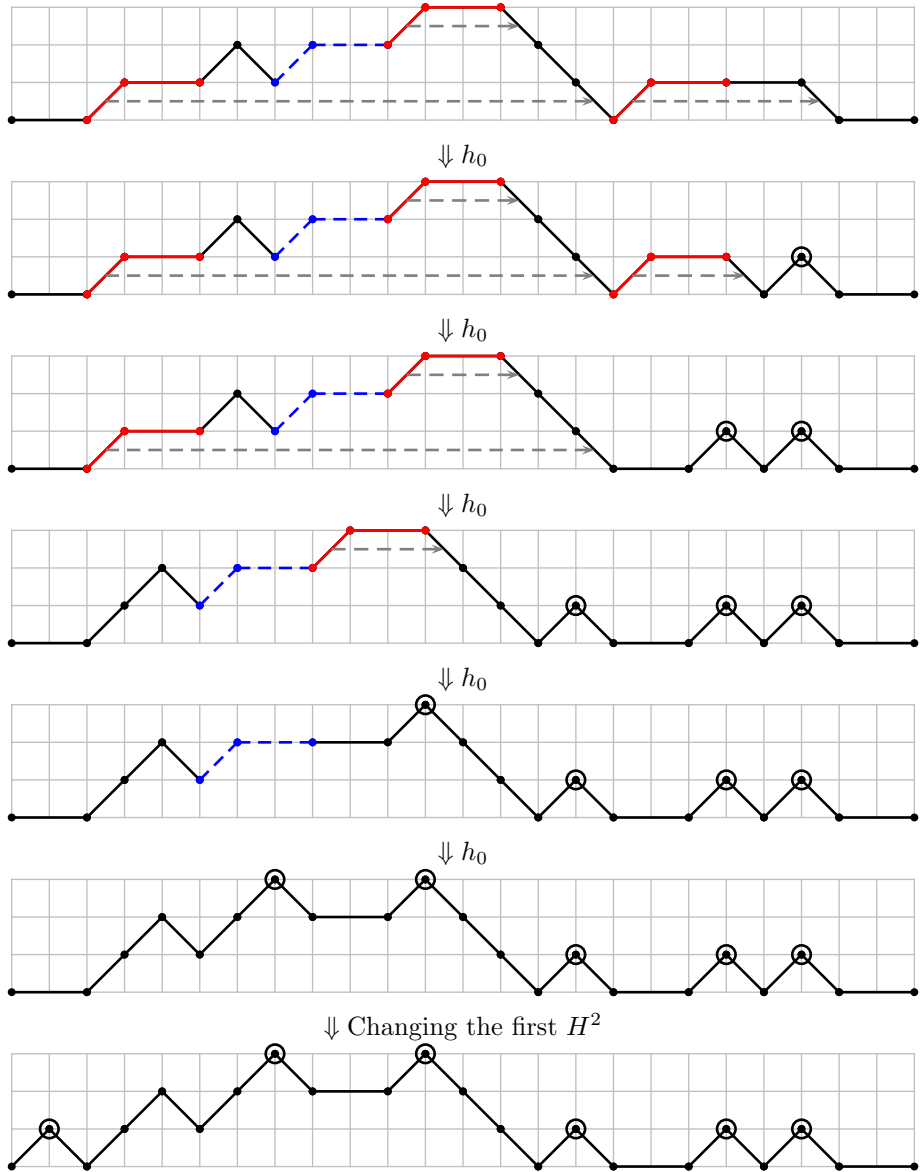

Now we define $h:\schodd{n}\rightarrow \schuhp{n}$ as follows.  Let
$L\in \schodd{n}$ and $L_0=L$. Then we define $L_i=h_0(L_{i-1})$ for
$i\geq1$ if $L_{i-1}$ is not UH-free. Since the number of UH-free
pairs of $L_i$ is one less than that of $L_{i-1}$, or they are the
same and
\[\xi(\mbox{the maximal UH-pair of } L_i) < \xi(\mbox{the maximal
  UH-pair of } L_{i-1}),\] we always get $L_r$ which is UH-free for
some $r$.  We define $h(L)$ to be $L_r$ if $L_r$ does not start with
$H^2$; and the lattice path obtained from $L_r$ by replacing $H^2$
with $UD$ otherwise.  For an example, see Figure~\ref{fig:h}.

\begin{thm}
The map $h:\schodd{n}\rightarrow \schuhp{n}$ is a bijection.
\end{thm}
\begin{proof}
In the procedure of $h$, the odd peaks are constructed from right to
left. Since $h_0$ is invertible, so is $h$.
\end{proof}

\subsection{The bijection $\iota:\schodd{n} \rightarrow \schevenp{n}$}

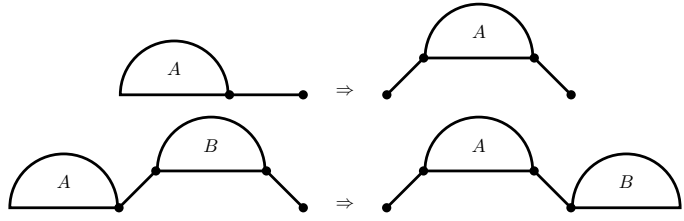
\begin{figure}
  \begin{center}
\scalebox{0.7}{\begin{pspicture}(0,0)(5,2)
\pswedge(1.5,0){1.5}{0}{180}
\rput(1.5,0.7){$A$}\DHS(3,0)
\end{pspicture}\hspace{0.5cm}
$\Rightarrow$\hspace{0.5cm}
\begin{pspicture}(0,0)(5,3)
\US(0,0) \pswedge(2.5,1){1.5}{0}{180}
\rput(2.5,1.7){$A$}\DS(4,1)
\end{pspicture}}

\scalebox{0.7}{\begin{pspicture}(0,0)(8,2)
\pswedge(1.5,0){1.5}{0}{180}
\rput(1.5,0.7){$A$}\US(3,0)
\pswedge(5.5,1){1.5}{0}{180}
\rput(5.5,1.7){$B$}\DS(7,1)
\end{pspicture}\hspace{0.5cm}
$\Rightarrow$\hspace{0.5cm}
\begin{pspicture}(0,0)(8,3)
\US(0,0) \pswedge(2.5,1){1.5}{0}{180}
\rput(2.5,1.7){$A$}\DS(4,1)
\pswedge(6.5,0){1.5}{0}{180} \rput(6.5,0.7){$B$}
\end{pspicture}}
  \end{center}
  \caption{The map $\iota$.}\label{fig:iota}
\end{figure}

For $L=S_1S_2\cdots S_k\in \schodd{n}$, we define $\iota(L)$ as follows.
\begin{enumerate}
\item If $S_k=H^2$, then
$\iota(L)=U S_1\cdots S_{k-1} D$. 
\item If $S_k=D$, then let $S_i$ be the up step corresponding to $S_k$
  and we define $\iota(L)=U S_1\cdots S_{i-1} D S_{i+1}\cdots S_{k-1}$.
\end{enumerate}
See Figure~\ref{fig:iota}.

Then $\iota(L)\in\schevenp{n}$. Clearly, $\iota:\schodd{n} \rightarrow
\schevenp{n}$ is a bijection.

\section{A direct bijection between $\NC{n}$ and $\avoiding{n}$}
\label{sec:direct}
Now we have a bijection $\phi \circ h \circ g \circ f \circ \psi:
\NCp{n}\rightarrow \avoidingp{n}$, see Figure~\ref{fig:bijections}. As noted in
the introduction, this induces a bijection between $\NC n$ and $\avoiding n$.
Since both $\NC{n}$ and $\avoiding{n}$ are partitions with some conditions, it
is natural to ask a direct bijection between them.  In this section we find such
a direct bijection.

From now on, we will identify a partition in $\avoiding{n}$ with its
canonical word.

A \emph{marked} partition is a partition in which each part may be marked.
Similarly a \emph{marked} word is a word in which each letter may be marked.

Let $\pi\in\NC{n}$.  For $i\in[n]$, let $T_i$ be the marked partition of $[i]$
obtained from $\pi$ by removing all the integers greater than $i$ and by marking
integers which are connected to an integer greater than $i$.  Using the sequence
$\emptyset=T_0, T_1, T_2,\ldots, T_n=\pi$ of marked partitions, we define a
sequence of marked words $\w_0, \w_1, \w_2,\ldots, \w_n$ as follows. Here, if
$(i,j)$ is an edge we say that \emph{$j$ is connected to $i$}. 

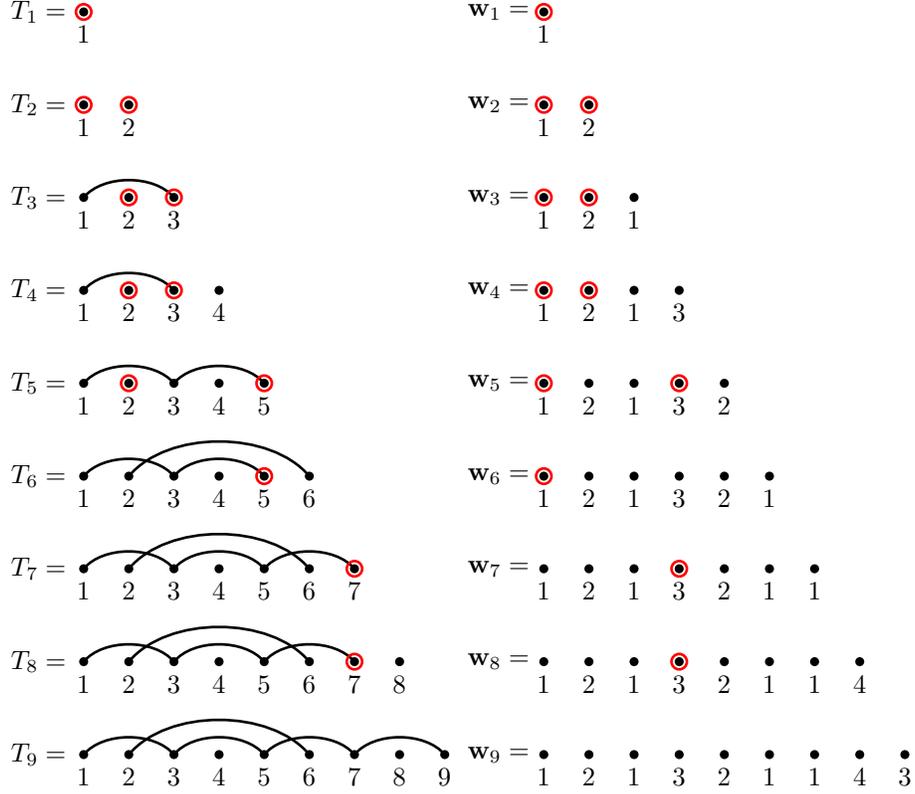
\begin{figure}
  \begin{center}
\psset{linewidth=1pt,unit=.6cm}
\begin{pspicture}(-1,.5)(9,2.5)
\rput(0,1){$T_1=$}
\vput{1}

\activevertex1
\end{pspicture}
\begin{pspicture}(-1,.5)(9,2.5)
\rput(0,1){$\w_1=$}
\cvput{1}[1]
\activevertex1
\end{pspicture}

\begin{pspicture}(-1,.5)(9,2.5)
\rput(0,1){$T_2=$}
\vput{1}\vput{2}

\activevertex1 \activevertex2
\end{pspicture}
\begin{pspicture}(-1,.5)(9,2.5)
\rput(0,1){$\w_2=$}
\cvput{1}[1]\cvput{2}[2]
\activevertex1 \activevertex2
\end{pspicture}

\begin{pspicture}(-1,.5)(9,2.5)
\rput(0,1){$T_3=$}
\vput{1}\vput{2}\vput{3}

\edge13 
\activevertex3 \activevertex2
\end{pspicture}
\begin{pspicture}(-1,.5)(9,2.5)
\rput(0,1){$\w_3=$}
\cvput{1}[1]\cvput{2}[2]\cvput{3}[1]
\activevertex1 \activevertex2
\end{pspicture}

\begin{pspicture}(-1,.5)(9,2.5)
\rput(0,1){$T_4=$}
\vput{1}\vput{2}\vput{3}\vput{4}

\edge13 
\activevertex3 \activevertex2
\end{pspicture}
\begin{pspicture}(-1,.5)(9,2.5)
\rput(0,1){$\w_4=$}
\cvput{1}[1]\cvput{2}[2]\cvput{3}[1]\cvput{4}[3]
\activevertex1 \activevertex2
\end{pspicture}

\begin{pspicture}(-1,.5)(9,2.5)
\rput(0,1){$T_5=$}
\vput{1}\vput{2}\vput{3}\vput{4}\vput{5}

\edge13 \edge35 
\activevertex5 \activevertex2
\end{pspicture}
\begin{pspicture}(-1,.5)(9,2.5)
\rput(0,1){$\w_5=$}
\cvput{1}[1]\cvput{2}[2]\cvput{3}[1]\cvput{4}[3]\cvput{5}[2]
\activevertex1 \activevertex4
\end{pspicture}

\begin{pspicture}(-1,.5)(9,2.5)
\rput(0,1){$T_6=$}
\vput{1}\vput{2}\vput{3}\vput{4}\vput{5}\vput{6}

\edge13 \edge35 \edge26 
\activevertex5
\end{pspicture}
\begin{pspicture}(-1,.5)(9,2.5)
\rput(0,1){$\w_6=$}
\cvput{1}[1]\cvput{2}[2]\cvput{3}[1]\cvput{4}[3]\cvput{5}[2]
\cvput{6}[1]
\activevertex1
\end{pspicture}

\begin{pspicture}(-1,.5)(9,2.5)
\rput(0,1){$T_7=$}
\vput{1}\vput{2}\vput{3}\vput{4}\vput{5}\vput{6}\vput{7}

\edge13 \edge35 \edge26 \edge57
\activevertex7
\end{pspicture}
\begin{pspicture}(-1,.5)(9,2.5)
\rput(0,1){$\w_7=$}
\cvput{1}[1]\cvput{2}[2]\cvput{3}[1]\cvput{4}[3]\cvput{5}[2]
\cvput{6}[1]\cvput{7}[1]
\activevertex4
\end{pspicture}

\begin{pspicture}(-1,.5)(9,2.5)
\rput(0,1){$T_8=$}
\vput{1}\vput{2}\vput{3}\vput{4}\vput{5}\vput{6}\vput{7}\vput{8}

\edge13 \edge35 \edge26 \edge57
\activevertex7
\end{pspicture}
\begin{pspicture}(-1,.5)(9,2.5)
\rput(0,1){$\w_8=$}
\cvput{1}[1]\cvput{2}[2]\cvput{3}[1]\cvput{4}[3]\cvput{5}[2]
\cvput{6}[1]\cvput{7}[1]\cvput{8}[4]
\activevertex4
\end{pspicture}

\begin{pspicture}(-1,.5)(9,2.5)
\rput(0,1){$T_9=$}
\vput{1}\vput{2}\vput{3}\vput{4}\vput{5}\vput{6}\vput{7}\vput{8}\vput{9} 

\edge13 \edge35 \edge26 \edge57 \edge79
\end{pspicture}
\begin{pspicture}(-1,.5)(9,2.5)
\rput(0,1){$\w_9=$}
\cvput{1}[1]\cvput{2}[2]\cvput{3}[1]\cvput{4}[3]\cvput{5}[2]
\cvput{6}[1]\cvput{7}[1]\cvput{8}[4]\cvput{9}[3]

\end{pspicture}

  \end{center}
  \caption{$T_i$'s and corresponding $\w_i$'s. Marked integers and
    marked letters are circled.}
\label{fig:direct}
\end{figure}

Let $\w_0$ be the empty word.  For $1\leq i\leq n$, $\w_i$ is defined as
follows.
\begin{enumerate}
\item If $i$ is not connected to any integer in $T_i$, then
  $\w_i=\w_{i-1}m$, where $m=\max(\w_{i-1})+1$. Otherwise, $i$ is
  connected to either the largest marked integer or the second largest
  marked integer of $T_{i-1}$.
  \begin{itemize}
  \item If $i$ is connected to the largest marked integer of
    $T_{i-1}$, then let $\w_i=\w_{i-1}a_1$, where $a_1$ is the
    rightmost marked letter of $\w_{i-1}$. And then we make the marked
    letter $a_1$ unmarked.
  \item If $i$ is connected to the second largest marked integer of
    $T_{i-1}$, then let $\w_i=\w_{i-1}a_2$, where $a_2$ is the second
    rightmost marked letter of $\w_{i-1}$. The second rightmost marked
    letter of $\w_{i-1}$ remains marked, however, we make the
    rightmost marked letter of $\w_{i-1}$ unmarked in $\w_i$.
  \end{itemize}
\item If $i$ is marked in $T_i$, then we find the largest letters in
  $\w_i$ and make the leftmost letter among them marked.
\end{enumerate}
For an example, see Figure~\ref{fig:direct}. 

\begin{lem}
  The word $\w_n$ obtained above is $12312$-avoiding.
\end{lem}
\begin{proof}
  Suppose $\w_n$ has a subsequence $abcab$ where $a<b<c$. When the second $b$ is
  added the first $b$ must have been marked. Moreover, the first $b$ must have
  been marked before adding the second $a$ because an unmarked integer becomes
  marked only if it is the largest integer (in this case at least $c$) in the
  sequence. Thus when the second $a$ is added, the first $a$ and $b$ have been
  marked. Since the first $a$ is the second rightmost marked integer at this
  moment, we must unmark the rightmost marked integer, the first $b$, and mark
  the largest integer which is at least $c$. Thus after this process, $b$ cannot
  be marked and we cannot have the second $b$, which is a contradiction.
\end{proof}

If we know $\w_n$, we can reverse this procedure.  For $1\leq i\leq n$,
$\w_{i-1}$ is obtained from $\w_i$ as follows. Suppose $m=\max(\w_i)$ and $t$ is
the last letter of $\w_i$.
\begin{enumerate}
\item If the leftmost $m$ is marked in $\w_i$, then make it unmarked.
\item If $t$ appears only once in $\w_i$ (equivalently $t$ is greater
  than any other letters in $\w_i$), then we simply remove
  $t$. Otherwise, find the leftmost $t$ in $\w_i$. 
  \begin{itemize}
  \item If the leftmost $t$ is unmarked, then we remove the last
    letter $t$ and make the leftmost $t$ marked.
  \item If the leftmost $t$ is marked, then we must have $t<m$ since
    we have made the leftmost $m$ unmarked. In this case we remove the
    last $t$, and make the leftmost $t$ still marked and the leftmost
    $m$ marked.
  \end{itemize}
\end{enumerate}

Now we construct $T_0,T_1,\ldots,T_n$ as follows.  Let $T_0=\emptyset$. For
$1\leq i\leq n$, $T_i$ is obtained as follows.
\begin{enumerate}
\item First, let $T_i$ be the marked partition obtained from $T_{i-1}$
  by adding $i$.
\item If the last letter of $\w_i$ is equal to the rightmost (resp. the second
  rightmost) marked letter of $\w_{i-1}$, then connect $i$ to the largest
  (resp. the second largest) marked integer, say $j$, of $T_{i-1}$, and make $j$
  unmarked.
\item Let $m=\max(\w_i)$. If the leftmost $m$ is marked in $\w_i$,
  then make $i$ marked in $T_i$.
\end{enumerate}

It is easy to check that this is the inverse map. Thus we get the
following theorem.

\begin{thm}
For $\pi\in\NC{n}$, the map $\pi\mapsto \w_n$ is a bijection from
$\NC{n}$ to $\avoiding{n}$.
\end{thm}

The bijection $\pi\mapsto \w_n$ is different from the composition $\phi \circ h
\circ g \circ f \circ \psi$. For instance, if $\pi=(\{1,3\},\{2\})$, then $\w_3=
121$ but $(\phi \circ h \circ g \circ f \circ \psi)(\pi)=112$.

Note that both $\NC n$ and $\avoiding n$ contain noncrossing partitions. It
would be interesting to find a bijection between $\NC n$ and $\avoiding n$
which sends noncrossings partitions to noncrossings partitions.


\end{document}